\documentclass{article}
\begin{document}
\newtheorem{proposition}{Proposition}[section]
\newtheorem{definition}{Definition}[section]
\newtheorem{lemma}{Lemma}[section]
\newcommand{\xl}{\stackrel{\rightharpoonup}{\cdot}}
\newcommand{\xr}{\stackrel{\leftharpoonup}{\cdot}}
\newcommand{\xlplus}{\stackrel{\rightharpoonup}{+}}
\newcommand{\xrplus}{\stackrel{\leftharpoonup}{+}}
\newcommand{\Ll}{\stackrel{\rightharpoonup}{L}}
\newcommand{\Lr}{\stackrel{\leftharpoonup}{L}}
\newcommand{\Rl}{\stackrel{\rightharpoonup}{R}}
\newcommand{\Rr}{\stackrel{\leftharpoonup}{R}}

\title{\bf Transformation Digroups}
\author{Keqin Liu\\Department of Mathematics\\The University of British Columbia\\Vancouver, BC\\
Canada, V6T 1Z2}
\date{September 11, 2004}
\maketitle

\begin{abstract} We introduce the notion of a transformation digroup and prove that 
every digroup is isomorphic to  a transformation digroup.
\end{abstract}

The purpose of this paper is to show how to choose a class of non-bijective 
transformations on a Cartesian product of two sets to define a transformation digroup 
on the Cartesian product. The main result of this paper is that every digroup is 
isomorphic to
a transformation digroup.

\medskip
The notion of a digroup we shall use in this paper was introduced in Chapter 6 of 
\cite{Liu2}. Its special case, which is the notion of a digroup with an identity, was 
introduced independently by people who work in different areas of mathematics 
(\cite{Felipe}, \cite{Kinyon1} and \cite{Liu1}). 

\medskip
After reviewing 
some basic definitions about digroups in Section 1, we introduce in Section 2 the 
notion of a symmetric digroup 
on the Cartesian product  $\Delta \times \Gamma$, where  $\Delta$ and  $\Gamma$ are 
two sets. 
If $|\Delta |=1$, then a 
symmetric digroup  on $\Delta \times \Gamma$ becomes the symmetric group on  the set 
$\Gamma$ , 
where $|\Delta |$ denotes the cardinality of $\Delta $. In Section 3 we prove that 
every digroup is isomorphic a subdigroup of a symmetric digroup, which is a better 
counterpart of Cayley's Theorem in the context of digroups. 

\section{The Notion of a Digroup}

The following definition of a digroup is a version of  Definition 6.1 of \cite{Liu2}.

\begin{definition}\label{def1.1} A nonempty set $G$ is called a {\bf digroup} if 
there are two binary operations $\xl$ and $\xr$ on $G$ such that the following three 
properties are satisfied.
\begin{description}
\item[(i)] ({\bf The Diassociative Law}) The two operations $\xl$ and $\xr$ are 
{\bf diassociative}; that is,
\begin{equation}\label{eq1.1} x\xl (y\xl z)=(x\xl y)\xl z=x\xl (y\xr z),\end{equation}
\begin{equation}\label{eq1.2} (x\xr y)\xl z=x\xr (y\xl z),\end{equation}
\begin{equation}\label{eq1.3} (x\xl y)\xr z=(x\xr y)\xr z=x\xr (y\xr z)\end{equation}
for all $x,y,z\in G$.
\item[(ii)] ({\bf Bar-unit}) There is an element $e$ of $G$ such that
\begin{equation}\label{eq1.4} x\xl e =x=e\xr x\end{equation}
for all $x\in G$.
\item[(iii)] ({\bf One-sided Inverses}) For each element $x$ in $G$, there exist two 
elements $x_{e}^{\stackrel{\ell}{-}1}$ and $x_{e}^{\stackrel{r}{-}1}$ of $G$ such that
\begin{equation}\label{eq1.5} 
x_e^{\stackrel{\ell}{-}1}\xl x=e=x\xr x_e^{\stackrel{r}{-}1}.
\end{equation}
\end{description}
\end{definition}

The diassociative law was introduced by J. -L. Loday to study Leibniz algebras. An 
element $e$ satisfying (\ref{eq1.4}) is called a {\bf bar-unit}, and the set of all 
bar-units is called the {\bf halo} (\cite{Loday2}). The two elements 
$x_e^{\stackrel{\ell}{-}1}$ and $x_e^{\stackrel{r}{-}1}$ are called a 
{\bf left inverse} and a {\bf right inverse} of $x$ with respect to the bar-unit $e$, 
respectively (\cite{Liu2}). 
The binary operations $\xl$ and $\xr$ are called the {\bf left product} and the 
{\bf right product}, respectively. We also use $(\,G, \,\xl , \,\xr \,)$ to signify 
that $G$ 
is a digroup with the left product $\xl$ and the right product $\xr$. 

\medskip
Both the left inverse $x_e^{\stackrel{\ell}{-}1}$ and the right inverse 
$x_e^{\stackrel{r}{-}1}$ of an element $x$ with respect to a bar-unit $e$ is unique 
(Proposition 6.1 in \cite{Liu2}), but $x_e^{\stackrel{\ell}{-}1}$ is generally not 
equal to $x_e^{\stackrel{r}{-}1}$ (Proposition 6.2 in \cite{Liu2}). 

\medskip
The next proposition shows that Definition~\ref{def1.1} of a digroup is independent of the choice of the bar-unit $e$.

\begin{proposition}\label{pr1.1} Let $G$ be a digroup. If $\alpha$ is a bar-unit of $G$, then every element $x$ of $G$ has both the left inverse $x_{\alpha}^{\stackrel{\ell}{-}1}$ and the right inverse $x_{\alpha}^{\stackrel{r}{-}1}$ with respect to $\alpha$; that is,
$$x_{\alpha}^{\stackrel{\ell}{-}1}\xl x=\alpha=x\xr x_{\alpha}^{\stackrel{r}{-}1}.$$

\end{proposition}

\medskip
\noindent
{\bf Proof} By Definition~\ref{def1.1}, there is a bar-unit $e$ of $G$ such that
$$x_e^{\stackrel{\ell}{-}1}\xl x=e=x\xr x_e^{\stackrel{r}{-}1}.$$

It follows that
$$(\alpha \xl x_e^{\stackrel{\ell}{-}1})\xl x=\alpha =x\xr (x_e^{\stackrel{r}{-}1}\xr \alpha),$$
which proves that
$$ x_{\alpha}^{\stackrel{\ell}{-}1}=\alpha \xl x_e^{\stackrel{\ell}{-}1} \quad
\mbox{and} \quad 
x_{\alpha}^{\stackrel{r}{-}1}=x_e^{\stackrel{r}{-}1}\xr \alpha.$$
\hfill\raisebox{1mm}{\framebox[2mm]{}}

\bigskip
If $G$ is a digroup, then the {\bf halo} of $G$ is denoted by $\hbar (G)$; that is
$$ \hbar (G): =\{\, \alpha \,\, |\,\, \mbox{$\alpha \in G$ and $x\xl \alpha =x=\alpha \xr x$ for all $x\in G$.} \, \}.$$

\begin{definition}\label{def1.0} Let $(\,G, \,\xl , \,\xr \,)$ be a digroup. An element $e$ of $G$ is called an {\bf identity} of $G$ if $e\in \hbar (G)$ and
$$e\xl x=x\xr e \quad\mbox{for all $x\in G$.}$$
\end{definition}

By Proposition 6.2 in \cite{Liu2}, if $e$ is a bar-unit of $G$, then $e$ is an 
identity of a digroup $G$ if and only if 
$x_e^{\stackrel{\ell}{-}1}=x_e^{\stackrel{r}{-}1}$ for all $x\in G$. Example 6 in  
\cite{Liu2} gives a digroup which does not have an identity.

\medskip
We now introduce the notion of a subdigroup.

\begin{definition}\label{def1.3} Let $H$ be a subset of a digroup $(\,G, \,\xl , \,\xr \,)$. If $H\cap \hbar (G)\neq\emptyset$ and $H$ is itself a digroup under the two binary operations of $G$, we say that $H$ is a {\bf subdigroup} of $G$. The notation $H \le _d G$ is used to indicate that $H$ is a subdigroup of $G$.
\end{definition}

The following definition introduces two important subdigroups of a digroup.

\begin{definition}\label{def1.4} Let $(\,G, \,\xl , \,\xr \,)$ be a digroup. 
\begin{description}
\item[(i)] The set 
$$\mathcal{Z}^t(G):=\{\, z\in G \, |\, \mbox{$z\xr x=x\xl z$ for all $x\in G$} \,\}$$
is called the {\bf target center} of $G$.
\item[(ii)] The set 
$$\mathcal{Z}^s(G):=\{\, z\in G \, |\, \mbox{$x\xr z=z\xl x$ for all $x\in G$} \,\}$$
is called the {\bf source center} of $G$.
\end{description}
\end{definition}

One can check that the target center $\mathcal{Z}^t(G)$ is a subdigroup of a digroup $G$ and
$$\hbar (G)\subseteq \mathcal{Z}^t (G). $$

If $G$ is a digroup with an identity, then the source center $\mathcal{Z}^s(G)$ is a subdigroup of $G$.

\medskip
We finish this section with the definition of an isomorphism between digroups.
 
\begin{definition}\label{def2.1} If $G$ and $\bar{G}$ are digroups, then a map $\varphi$ from $G$ to $\bar{G}$ is called an {\bf isomorphism} if $\varphi$ is bijective and 
$$\varphi (x\ast y)=\varphi (x)\ast\varphi (y)$$
for all $x,y\in G$ and $\ast =\xl , \; \xr$.
\end{definition}

\bigskip
\section{Symmetric Digroups}

Let $\mathcal{T} (\Delta\times \Gamma)$ be the set of all maps from $\Delta\times \Gamma$ to 
$\Delta\times \Gamma$, where $\Delta$ and $\Gamma$  are two sets. Then $\mathcal{T} (\Delta\times \Gamma)$ is a semigroup with the identity $1$ under the product:
$$ fg:=f\cdot g, $$
where $1$ is the identity map, and the product $f\cdot g$ is the composite of $g$ and $f$ ( $f$ following $g$ ):
$$ (f\cdot g)(x):=f(g(x)) \quad\mbox{for $x\in \Delta\times \Gamma$}.$$

\begin{definition}\label{def4.1} An element $\ell$ of $\mathcal{T} (\Delta\times \Gamma)$ is called a {\bf $\ell$-map} on $\Delta\times \Gamma$ if there exists  an element $(s, f)\in \Delta\times Sym\Gamma$ such that 
$$
\ell (k, i)=(s, f(i)) \quad\mbox{for $(k,i)\in \Delta\times \Gamma$},
$$
where $Sym\Gamma$ is the symmetric group on $\Gamma$. 
\end{definition}

It is clear that the element $(s, f)$ in Definition~\ref{def4.1} is determined uniquely by the $\ell$-map $\ell$. Hence, $\ell$ in Definition~\ref{def4.1} is also called the $\ell$-map on $\Delta\times \Gamma$ induced by $(s, f)\in \Delta\times Sym\Gamma$. We shall use the notation 
$\ell _{s, f}$ to indicate that the $\ell$-map induced by $(s, f)$. Thus, we have
\begin{equation}\label{eq4.1}
\ell _{s, f} (k, i)=(s, f(i)) \quad\mbox{for $(k,i)\in \Delta\times \Gamma$}.
\end{equation}

\medskip
Let $\ell _{s, f}$ and $\ell _{t, g}$ be two $\ell$-maps on $\Delta\times \Gamma$. If $(k, i)\in \Delta\times \Gamma$, then
$$\left(\ell _{s, f}\ell _{t, g}\right)(k, i)
=\ell _{s, f}(t, g(i))=(s, fg(i))=\ell _{s, fg}(k, i)$$
by (\ref{eq4.1}). Hence, we have
\begin{equation}\label{eq4.2}
\ell _{s, f}\ell _{t, g}=\ell _{s, fg} \quad\mbox{for $(s,f)$, $(t,g)\in \Delta\times Sym\Gamma$}.
\end{equation}

\begin{proposition}\label{pr4.1} If $\mathcal{G}$ is a subgroup of $Sym\Gamma$, then the set
$$ \ell _{\Delta\times \mathcal{G}}:=\{\,\ell _{s, f} \, | \, (s,f)\in \Delta\times \mathcal{G} \,\} $$ 
is a subsemigroup of $\mathcal{T} (\Delta\times Sym\Gamma)$ having the following two properties:
\begin{enumerate}
\item $\ell _{\Delta\times \mathcal{G}}$ has a right unit $e:=\ell _{0,1}$, where $0$ is a fixed element of $\Delta$.
\item Every element $\ell _{s, f}$ of $\ell _{\Delta\times \mathcal{G}}$ has a left inverse $\ell _{s,f}^{\stackrel{\ell}{-}1}:=\ell _{0, f^{-1}}$ in 
$\ell _{\Delta\times \mathcal{G}}$ with respect to the right unit $e=\ell _{0,1}$.
\end{enumerate}
\end{proposition}

\medskip
\noindent
{\bf Proof} It is clear by (\ref{eq4.2}).
\hfill\raisebox{1mm}{\framebox[2mm]{}}

\begin{proposition}\label{pr4.2} Let $\mathcal{G}$ be a subgroup of $Sym\Gamma$, and 
let $\theta : \mathcal{G} \to  Sym\Delta$ be a group homomorphism. 
\begin{description}
\item[(i)] The map
$$ f\mapsto \bar{\ell} _f \quad\mbox{ for $f\in \mathcal{G}$} $$
is a group homomorphism from $\mathcal{G}$ to $Sym(\Delta \times \Gamma)$, where $\bar{\ell} _f : \Delta \times \Gamma \to \Delta \times \Gamma $ is defined by
\begin{equation}\label{eq4.3}
\bar{\ell} _f (k, i):=\left( \theta (f)(k), f(i) \right) \quad\mbox{for $(k, i)\in \Delta \times \Gamma$ }.
\end{equation}
$\bar{\ell} _f$ will be called a {\bf $\theta$-permutation} on $\Delta \times \Gamma$.
\item[(ii)] If $(s, f)$, $(t, g) \in \Delta\times \Gamma$, then
\begin{eqnarray}\label{eq4.4}
\bar{\ell} _f \ell _{t, g}&=&\ell _{\theta (f)(t), fg}, \\&&\nonumber\\
\label{eq4.5}
\ell _{s, f} \bar{\ell} _g&=&\ell _{s, fg} .
\end{eqnarray}
\end{description}
\end{proposition}

\medskip
\noindent
{\bf Proof} (i) It is clear that $\bar{\ell} _f\in Sym(\Delta \times \Gamma)$ by (\ref{eq4.3}). For $f$, $g\in\mathcal{G}$ and $(k, i)\in \Delta \times \Gamma$, we have
\begin{eqnarray*}
&&\left( \bar{\ell} _f \bar{\ell} _g\right)(k,i)\\
&&\\&=&\bar{\ell} _f\left( \theta (g)(k), g(i) \right)=
\left( \theta (f)\theta (g)(k), fg(i)\right)\\
&&\\&=&\left( \theta (fg)(k), (fg)(i)\right)=\bar{\ell} _{fg}(k,i),
\end{eqnarray*}
which implies that
\begin{equation}\label{eq4.6}
\bar{\ell} _f \bar{\ell} _g=\bar{\ell} _{fg} \quad\mbox{for $f$, $g\in\mathcal{G}$}.
\end{equation}
This proves that the map $ f\mapsto \bar{\ell} _f$ is a group homomorphism. 

\medskip
(ii) For $(k, i)\in \Delta \times \Gamma$, we have
$$
\left( \bar{\ell} _f \ell _{t, g}\right) (k, i)=\bar{\ell} _f (t, g(i))=
\left( \theta (f)(t), (fg)(i) \right)=\ell _{\theta (f)(t), fg} (k, i)
$$
and
$$
\left( \ell _{s, f} \bar{\ell} _g \right) (k, i)=
\ell _{s, f} \left( \theta (g)(k), g(i) \right)=
\left( s, (fg)(i) \right)=\ell _{s, fg} (k, i).
$$
Hence, (ii) is true.
\hfill\raisebox{1mm}{\framebox[2mm]{}}

\bigskip
Motivated by the facts in \cite{Liu1}, we introduce the construction of a transformation digroup in the next proposition.

\begin{proposition}\label{pr4.3} If $\mathcal{G}$ is a subgroup of $Sym\Gamma$ and $\theta : \mathcal{G} \to  Sym\Delta$  is a group homomorphism, then
 $\ell _{\Delta\times \mathcal{G}}$ is a digroup under the following two binary operations:
\begin{eqnarray}
\label{eq4.7}&&\ell _{s, f} \xl \ell _{t, g} :=\ell _{s, f} \ell _{t, g}=\ell _{s, fg},\\
\label{eq4.8}&&\ell _{s, f} \xr \ell _{t, g} :=
\bar{\ell} _f \ell _{t, g}=\ell _{\theta (f)(t), fg},
\end{eqnarray}
where $(s, f)$, $(t, g)\in \Delta\times \mathcal{G}$. $e:=\ell _{0, 1}$ is a bar-unit, and the left inverse and the right inverse of an element with respect to the bar-unit  $\ell _{0, 1}$ are given by
\begin{equation}\label{eq4.9}
\left(\ell _{s,f}\right)_e^{\stackrel{\ell}{-}1}=\ell _{0, f^{-1}}, \qquad
\left(\ell _{s,f}\right)_e^{\stackrel{r}{-}1}=\ell _{\theta (f^{-1})(0), f^{-1}}.
\end{equation}
\end{proposition}

\medskip
\noindent
{\bf Proof} For $(s, f)$, $(t, g)$ and $(v, h) \in \Delta\times \mathcal{G}$, we have
$$
\ell _{s, f} \xl \left(\ell _{t, g}\xl \ell _{v, h}\right)=
 \left(\ell _{s, f} \xl \ell _{t, g}\right)\xl \ell _{v, h}=
\ell _{s, f} \xl \left(\ell _{t, g}\xr \ell _{v, h}\right)=\ell _{s, fgh},
$$
$$
\left(\ell _{s, f} \xr \ell _{t, g}\right)\xl \ell _{v, h}=
\ell _{s, f} \xr \left(\ell _{t, g}\xl \ell _{v, h}\right)=\ell _{\theta (f)(t), fgh},
$$
$$
\left(\ell _{s, f} \xl \ell _{t, g}\right)\xr \ell _{v, h}=
 \left(\ell _{s, f} \xr \ell _{t, g}\right)\xr \ell _{v, h}=
\ell _{s, f} \xr \left(\ell _{t, g}\xr \ell _{v, h}\right)=\ell _{\theta (fg)(v), fgh}.
$$

This proves that the two binary operations $\xl$ and $\xr$ are diassociative. The remaining parts are clear by (\ref{eq4.7}) and (\ref{eq4.8}).
\hfill\raisebox{1mm}{\framebox[2mm]{}}

\medskip
A subdigroup of $\ell _{\Delta\times \mathcal{G}}$ is called the {\bf transformation digroup} on $\Delta \times \Gamma$ induced by $\left( \mathcal{G}, \theta \right)$. In particular, the digroup $\ell _{\Delta\times Sym\Gamma}$ is called the {\bf symmetric digroup} on $\Delta \times \Gamma$ induced by the group homomorphism $\theta :  Sym\Gamma \to  Sym\Delta $. It is clear that 
the symmetric digroup becomes the symmetric group on $\Gamma$ if $|\Delta |=1$.

\medskip
We finish this section with the description of the halo and subdigroups of $\ell _{\Delta\times \mathcal{G}}$.

\begin{proposition}\label{pr4.4} Let $\ell _{\Delta\times \mathcal{G}}$ be the transformation digroup on $\Delta \times \Gamma$ induced by $\left( \mathcal{G}, \theta \right)$. Let $\mathcal{Z} (\mathcal{G})$ be the center of the group $\mathcal{G}$. Then
\begin{description}
\item[(i)] $\hbar \left( \ell _{\Delta\times \mathcal{G}} \right)=
\{ \, \ell _{s, 1} \, | \, s\in\Delta \, \}$.
\item[(ii)] $\ell _{s, 1}$ is an identity of $\ell _{\Delta\times \mathcal{G}}$ if and only if 
$Im\theta \subseteq \left(Sym\Delta\right)_s$, where
$$ \left(Sym\Delta\right)_s: =\{\, x\in Sym\Delta \, | \, x(s)=s \, \} $$
is the {\bf stabilizer} of $s\in\Delta$ in $Sym\Delta$.
\item[(iii)] $H$ is a subdigroup of $\ell _{\Delta\times \mathcal{G}}$ if and only if there exist a subgroup $\mathcal{H}$ of $\mathcal{G}$ and a fixed block $\Omega$ of $Im \mathcal{H}$ such that $H=\ell _{\Omega\times \mathcal{H}}$.
\item[(iv)] The target center $\mathcal{Z}^t \left(\ell _{\Delta\times \mathcal{G}}\right)$ is given by
$$\mathcal{Z}^t \left(\ell _{\Delta\times \mathcal{G}}\right):=\{\, \ell _{s, f} \, |\, \mbox{$s\in \Delta$ and $f\in Ker \theta \cap \mathcal{Z} (\mathcal{G})$} \,\}.$$
\item[(v)] The source center $\mathcal{Z}^s \left(\ell _{\Delta\times \mathcal{G}}\right)$ is given by
$$\mathcal{Z}^s \left(\ell _{\Delta\times \mathcal{G}}\right):=\{\, \ell _{s, f} \, |\, 
\mbox{$s\in \Delta$, $Im \theta\subseteq \left(Sym\Delta\right)_s$ and 
$f\in \mathcal{Z} (\mathcal{G})$} \,\}.$$
\end{description}
\end{proposition}

\medskip
\noindent
{\bf Proof} The results follow from (\ref{eq4.7}) and (\ref{eq4.8}).
\hfill\raisebox{1mm}{\framebox[2mm]{}}

\section{The Counterpart of Cayley's Theorem}

In this section we prove that every digroup is isomorphic to a transformation digroup. We begin with the following property of left translations\footnote{The basic properties of left translations were given in Chapter 2 of \cite{Liu2}}.

\begin{proposition}\label{pr4.5} Let $e$ be a bar-unit of a digroup $G$ and let
$\phi : G\to \stackrel{\rightharpoonup}{e}(G)$ be a map defined by
$$\phi (g): =\vec{g},$$
where $g\in G$, $\vec{g}:=e\xl g$ and
$$
\stackrel{\rightharpoonup}{e}(G): =\left\{\, \vec{g}_i \, \left| \, \mbox{$i\in \bar{\Gamma}$ and $\vec{g}_i\ne \vec{g}_j$ for $i$, $j\in \bar{\Gamma}$ and $i\ne j$}\right. \, \right\}.
$$
If $f\in G$, then $\stackrel{\leftharpoonup}{L}_f\in Sym (\Gamma)$, where
$$
\Gamma: =\left\{\, \left. \phi ^{-1}\left(\vec{g}_i\right) \, \right | \, i\in \bar{\Gamma} \, \right\}.
$$
\end{proposition}

\medskip
\noindent
{\bf Proof} It is clear that $\phi$ is surjective. Hence, $\phi ^{-1}\left(\vec{g}_i\right)$ is not empty for all $i\in  \bar{\Gamma}$. Since
\begin{eqnarray*}
&& e\xl \stackrel{\leftharpoonup}{L}_f\left( \phi ^{-1}\left(\vec{g}_i\right)\right)=
e\xl \left(f\xr \phi ^{-1}\left(\vec{g}_i\right) \right)\\
&=&\left(e\xl f \right)\xl \left(e\xl \phi ^{-1}\left(\vec{g}_i\right)\right)=
e\xl f\xl \vec{g}_i=\vec{f}\xl \vec{g}_i,
\end{eqnarray*}
we have
\begin{equation}\label{eq4.10}
\stackrel{\leftharpoonup}{L}_f\left( \phi ^{-1}\left(\vec{g}_i\right)\right)
\subseteq \phi ^{-1} \left( \vec{f}\xl \vec{g}_i\right).
\end{equation}

Using the properties of left translations and (\ref{eq4.10}), we have
\begin{eqnarray*}
&&\phi ^{-1} \left( \vec{f}\xl \vec{g}_i\right)=
\stackrel{\leftharpoonup}{L}_e\left( \phi ^{-1} \left( \vec{f}\xl \vec{g}_i\right)\right)=
\stackrel{\leftharpoonup}{L}_{f\xr f_{e}^{\stackrel{r}{-}1}}\left( \phi ^{-1} \left( \vec{f}\xl \vec{g}_i\right)\right)\\
&=&\stackrel{\leftharpoonup}{L}_{f}\stackrel{\leftharpoonup}{L}_{f_{e}^{\stackrel{r}{-}1}}
\left( \phi ^{-1} \left( \vec{f}\xl \vec{g}_i\right)\right)\subseteq
\stackrel{\leftharpoonup}{L}_{f}
\left(\phi ^{-1} \left( e\xl f_{e}^{\stackrel{r}{-}1}\xl \vec{f}\xl \vec{g}_i\right) \right)\\
&=&\stackrel{\leftharpoonup}{L}_{f}
\left(\phi ^{-1} \left( f_{e}^{\stackrel{\ell}{-}1}\xl \vec{f}\xl \vec{g}_i\right) \right)=
\stackrel{\leftharpoonup}{L}_{f}\left(\phi ^{-1} \left( \vec{g}_i\right) \right)
\end{eqnarray*}
or
\begin{equation}\label{eq4.11}
\stackrel{\leftharpoonup}{L}_f\left( \phi ^{-1}\left(\vec{g}_i\right)\right)
\supseteq \phi ^{-1} \left( \vec{f}\xl \vec{g}_i\right).
\end{equation}

It follow from (\ref{eq4.10}) and (\ref{eq4.11}) that
\begin{equation}\label{eq4.12}
\stackrel{\leftharpoonup}{L}_f\left( \phi ^{-1}\left(\vec{g}_i\right)\right)
=\phi ^{-1} \left( \vec{f}\xl \vec{g}_i\right) \quad\mbox{for $i\in \bar{\Gamma}$}.
\end{equation}

Since $\stackrel{\leftharpoonup}{L}_f$ is a permutation on $G$, $\stackrel{\leftharpoonup}{L}_f$ can be regarded as a permutation on $\Gamma$ by (\ref{eq4.12}).
\hfill\raisebox{1mm}{\framebox[2mm]{}}

\medskip
For an element $f$ of a digroup $G$, we define a map $\Psi _f : G\to G$ by
\begin{equation}\label{eq3.1}
\Psi _f (x): =f\xr x\xl f_e^{\stackrel{\ell}{-}1}=f\xr x\xl f_e^{\stackrel{r}{-}1}
\quad\mbox{for $ x\in G$,}
\end{equation}
where $e$ is a bar-unit of $G$. Since the definition of  $\Psi _f$ is independent of 
the choice of the bar-unit $e$, (\ref{eq3.1}) is also written as
$$\Psi _f (x): =f\xr x\xl f ^{\stackrel{\ell}{-}1}=f\xr x\xl f^{\stackrel{r}{-}1}
\quad\mbox{for $x\in G$,}$$
where $f^{\stackrel{\ell}{-}1}$ and $f^{\stackrel{r}{-}1}$ denote the left inverse 
and the right inverse of $f$ with respect to any bar-unit, respectively.

\medskip
It is clear that $\Psi _f$ can be regarded as a permutation on the halo of $G$.

\begin{proposition}\label{pr4.6} If $G$ is a digroup, then
\begin{equation}\label{eq4.13}
\theta : \stackrel{\leftharpoonup}{L}_f \mapsto \left.\Psi _f\right| \hbar (G) \quad
\mbox{for $f\in G$}
\end{equation}
is a group homomorphism from the subgroup $\stackrel{\leftharpoonup}{L}_G$ of 
$Sym \Gamma$ to $Sym \,\hbar (G)$, where
$$ \stackrel{\leftharpoonup}{L}_G:=\{ \, \stackrel{\leftharpoonup}{L}_f \, |\, f\in G \, \} \quad\mbox{and}\quad 
\Gamma: =\left\{\, \left. \phi ^{-1}\left(\vec{g}_i\right) \, \right | \, i\in \bar{\Gamma} \, \right\}$$
\end{proposition}

\medskip
\noindent
{\bf Proof} $\theta$ is well-defined. In fact, if $f$, $g\in G$ and $\stackrel{\leftharpoonup}{L}_f=\stackrel{\leftharpoonup}{L}_g$, then $f\xr e= g\xr e$ for a 
bar-unit $e$ of $G$.  For all $x\in G$, we have
\begin{eqnarray*}
&&\Psi _f(x)=f\xr x\xl f_e^{\stackrel{\ell}{-}1}\\
&=&f\xr e\xr x\xl e\xl f_e^{\stackrel{\ell}{-}1}=
\left(f\xr e\right) \xr x\xl \left(f\xr e\right)_e^{\stackrel{\ell}{-}1}\\
&=&\left(g\xr e\right) \xr x\xl \left(g\xr e\right)_e^{\stackrel{\ell}{-}1}=\Psi _g(x).
\end{eqnarray*}
Hence, $\Psi _f=\Psi _g$.

\medskip
Since 
$$ \Psi _{fg}=\Psi _f\Psi _g \quad\mbox{for $f$, $g\in G$},$$
$\theta$ is a group homomorphism.
\hfill\raisebox{1mm}{\framebox[2mm]{}}

\bigskip
Let $G$ be a digroup. By Proposition~\ref{pr4.5} and Proposition~\ref{pr4.6}, 
$\ell _{\Delta\times \stackrel{\leftharpoonup}{L}_G}$ is a transformation digroup on $\Delta\times \Gamma$ induced by $\left( \mathcal{G}, \theta \right)$, where 
$$\Delta:=\hbar (G) \quad \mbox{and}\quad
\Gamma: =\left\{\, \left. \phi ^{-1}\left(\vec{g}_i\right) \, \right | \, i\in 
\bar{\Gamma} \, \right\}.$$ 

We now define a map $\lambda$ from $G$ to $\ell _{\Delta\times \stackrel{\leftharpoonup}{L}_G}$ by
\begin{equation}\label{eq4.14}
\lambda : \alpha \xl f \mapsto \ell _{\alpha , \stackrel{\leftharpoonup}{L}_f }\quad\mbox{for $\alpha \in \hbar (G)$ and $f\in G$.}
\end{equation}

\medskip
First, we prove that $\lambda$ is well-defined. Every element of $G$ is of the form $\alpha \xl f$ for some $\alpha \in \hbar (G)$ and $f\in G$. If $\alpha \xl f=\beta \xl g$, where $\alpha,\beta \in \hbar (G)$ and $f$, $g\in G$, then $\alpha =\beta$ and 
$\left( \alpha \xl f \right) \xr \alpha =\left( \alpha \xl g\right) \xr \alpha $ or 
$f\xr \alpha=g\xr \alpha$. It follows that $\stackrel{\leftharpoonup}{L}_f=
\stackrel{\leftharpoonup}{L}_g$. This proves that $\lambda$ is well-defined.

\medskip
Next, we prove that $\lambda$ is injective. For $\alpha$, $\beta \in\hbar (G)$ and $f$, $g\in G$, we have
\begin{eqnarray*}
&&\lambda \left( \alpha \xl f \right)=\lambda \left( \beta \xl g \right)\\
&\Rightarrow& 
\ell _{\alpha , \stackrel{\leftharpoonup}{L}_f }=\ell _{\beta , \stackrel{\leftharpoonup}{L}_g}\\
&\Rightarrow& \mbox{$\alpha =\beta$ and $\stackrel{\leftharpoonup}{L}_f=\stackrel{\leftharpoonup}{L}_g$}\\
&\Rightarrow& \mbox{$\alpha =\beta$ and $f\xr \alpha =g\xr \alpha $}\\
&\Rightarrow& \mbox{$\alpha =\beta$ and $\alpha \xl \left( f\xr \alpha\right) =\alpha \xl \left( g\xr \alpha \right)$}\\
&\Rightarrow& \alpha \xl f =\beta \xl g.
\end{eqnarray*}

\medskip
Finally, we prove that $\lambda$ preserves both the left product and the right product on $G$. For $\alpha$, $\beta \in\hbar (G)$ and $f$, $g\in G$, we have
\begin{eqnarray*}
&& \lambda \left(\left( \alpha \xl f \right) \xl \left(\beta \xl g \right)\right)\\&&\\
&=&\lambda \left(\alpha \xl f \xl g \right)=
\ell _{\alpha , \stackrel{\leftharpoonup}{L}_{f\xl g}}\\&&\\
&=&\ell _{\alpha , \stackrel{\leftharpoonup}{L}_{f}\stackrel{\leftharpoonup}{L}_{g}}
=\ell _{\alpha , \stackrel{\leftharpoonup}{L}_{f}}\xl
\ell _{\beta , \stackrel{\leftharpoonup}{L}_{g}}\\&&\\
&=&\lambda \left( \alpha \xl f \right)\xl \lambda \left(\beta \xl g \right)
\end{eqnarray*}
and
\begin{eqnarray*}
&& \lambda \left(\left( \alpha \xl f \right) \xr \left(\beta \xl g \right)\right)\\&&\\
&=&\lambda \left(\alpha \xr f\xr \beta  \xl g \right)\\&&\\
&=&\lambda \left(\left(f\xr \beta \xl f^{\stackrel{\ell}{-}1}\right)\xl f \xl g \right)\\&&\\
&=&\ell _{ f\xr \beta \xl f^{\stackrel{\ell}{-}1},\,\, \stackrel{\leftharpoonup}{L}_{f\xl g}}=
\ell _{ \Psi _f \left(\beta\right), \,\, \stackrel{\leftharpoonup}{L}_{f}\stackrel{\leftharpoonup}{L}_{g}}\\&&\\
&=&\ell _{ \theta \left(\stackrel{\leftharpoonup}{L}_{f}\right) \left(\beta\right), \,\, \stackrel{\leftharpoonup}{L}_{f}\stackrel{\leftharpoonup}{L}_{g}}=
\ell _{\alpha , \stackrel{\leftharpoonup}{L}_{f}}\xr
\ell _{\beta , \stackrel{\leftharpoonup}{L}_{g}}\\&&\\
&=&\lambda \left( \alpha \xl f \right)\xr \lambda \left(\beta \xl g \right).
\end{eqnarray*}

Thus we get the following counterpart of Cayley's Theorem.

\begin{proposition}\label{pr4.7} Any digroup is isomorphic to a transformation digroup.
\hfill\raisebox{1mm}{\framebox[2mm]{}}
\end{proposition}

\bigskip
\bigskip
{\bf Acknowledgment} I would like to thank Michael K. Kinyon for telling me  
the 
references 
\cite{Felipe}, \cite{Kinyon1} and \cite{Kinyon} after I published \cite{Liu2}.


\bigskip
\begin{thebibliography}{99}
\bibitem{Felipe} R. Felipe, \textsl{Generalized Loday algebras and digroups}, communication Tecnica No I-04-01/21-01-2004. www.cimat.mx/reportes/enlinea/I-04-01.pdf
\bibitem{Hall} Marshall Hall, Jr, \textsl{The theory of groups}, The Macmillan Company, 1959.
\bibitem{Kinyon1} Michael K. Kinyon, \textsl{The coquecigrue of a Leibniz algebra}, presented at {\it Alan Fest}, a conference in honor of the 60th birthday of Alan Weinstein, Erwin Schr$\ddot{o}$dinger Institute, Vienna, Austria, 4 August 2003, www.impa.br/jair/alanposter/coquecigrue.pdf
\bibitem{Kinyon} Michael K. Kinyon, \textsl{Leibniz Algebras, Lie Racks, and digroups}, arXiv: math. RA/0403509v2 31 Mar 2004.
\bibitem{Liu1} Keqin Liu, \textsl{ A class of group-like objects}, arXiv: math. RA/0311396v1 22 Nov 2003.
\bibitem{Liu2} Keqin Liu, \textsl{The generalizations of groups}, 153 Publishing, 2004.
\bibitem{Loday2} J.-L.Loday, A.Frabetti, F.Chapoton, F.Goichot, \textsl{Dialgebras and related operads}, Lecture Notes in Mathematics 1763, Springer, 2000

\end{thebibliography}
\end{document}